\documentclass[12pt]{scrartcl}

\pdfoutput=1
\usepackage[utf8]{inputenc} 
\usepackage{mathtools} 
\usepackage{amssymb}

\usepackage[square, numbers]{natbib}

\usepackage{amsfonts,amsmath,amsthm,amscd} 
\usepackage[british]{babel} 

\newtheorem{theorem}{Theorem}

\title{Cocycles in Lie Groups, Cochains and Regularity Problem} 
\author{\LARGE{Rosário D. Laureano}}

\bigskip

\date {ISCTE - Instituto Universitário de Lisboa \\ ISTAR - Information Sciences, Technologies and Architecture Research Centre }

\usepackage{graphicx}

\begin{document}

\maketitle

%\bigskip

%Project Title: \textbf{Research on Cohomological Dynamics and Applications in Quantum Complex Networks}

\textbf{ABSTRACT:} After the fundamental work of Livschitz in \cite{L1,L2}, various research directions emerged, among which the following stand out:
(i) the study of cocycles with values in groups and semigroups beyond $\mathbb{R}$, as well as the investigation of corresponding regularity results;
(ii) the analysis of how a certain degree of regularity ($C^{k}$ for $k=1,2,\ldots,\infty ,\omega$) of the cocycle can confer corresponding regularity to the solution of the cohomological equation; and
(iii) the study of higher-dimensional cohomology naturally associated with the action of groups such as $\mathbb{Z}^{k}$ or $\mathbb{R}^{k}$.
The aim of this article is to present, as self-contained as possible, a review of the natural generalizations of the notions of cocycles and cochains, as well as their corresponding results, in the study of cohomological equations.

\section{Introduction}

\qquad In dynamical systems theory, various problems of considerable importance can be reduced to solving an equation of the form: \begin{equation} \label{um} \varphi=\Phi\circ f -\Phi \end{equation} where $f:X\to X$ is a dynamical system, $\varphi :X\to \mathbb{R}$ is a known function and $\Phi:X\to \mathbb{R}$ is unknown. Equation \eqref{um} is called a \emph{cohomological equation}.

The study of cohomological equations is related, in particular, to the analysis of conjugacies to irrational rotations of the circle, the existence of absolutely continuous measures for expanding transformations of the circle, and the topological stability of hyperbolic automor-phisms of the torus. Such equations also naturally arise in statistical mechanics and celestial mechanics. 

Some results established by Livschitz in the 1970s \cite{L2} address precisely the possibility of obtaining solutions to cohomological equations in the context of hyperbolic dynamics. For a hyperbolic dynamical system, Livschitz's theorem provides a necessary and sufficient condition, based solely on information from periodic orbits, for the existence of H\"{o}lder solutions. This is one of the main tools for obtaining global cohomological information from periodic data.

A pragmatic approach to Livschitz's theorem, oriented towards the study of cohomology in dynamical systems, was presented in \cite{RL1,RL2}, emphasizing the relationship between the existence of solutions to cohomological equations and the behaviour of cocycles along periodic orbits. Following a preliminary demonstration of Anosov's closing lemma for hyperbolic diffeomorphisms, a detailed proof of Livschitz's theorem for hyperbolic diffeomorphisms is provided in \cite{RL2}, closely following the approach of Katok and Hasselblatt in \cite{KH}. The only published proof of Livschitz's theorem for flows is by Livschitz himself \cite{L2}. In \cite{RL1}, a proof of Livschitz's theorem for the continuous-time case is given, and its generalisation to suspension flows is discussed; this generalisation enabled a second proof of Livschitz's theorem for flows, based on the construction of Markov partitions by Bowen \cite{Bo} and Ratner \cite{R} for hyperbolic flows.

In this article, we present a review of cohomology in dynamical systems, with the aim of understanding the concept of cocycles and coboundaries more generally, in groups beyond $\mathbb{Z}$ or $\mathbb{R}$, and in higher dimensions. The problem of regularity of the solution of a cohomological equation, stemming from the regularity of the cocycle, is addressed, and regularity results are presented.

The study of cohomology in dynamical systems becomes more complex in non-hyperbolic dynamics, even when dealing with the cohomology of actions of $\mathbb{Z}$ or $\mathbb{R}$. Veech obtained an important result in this direction in \cite{V}, where he establishes the $C^{\infty }$ Livschitz's property for partially hyperbolic endomorphisms of the torus (i.e. showing that $C^{\infty }$ cocycles satisfying conditions relative to periodic orbits possess $C^{\infty }$ trivializations). This research focused solely on the hyperbolic case, presenting a selection of results that facilitate entry into theory.

\section{Cocycles with Values in Lie Groups}

\qquad Let $G$ be a group that acts on a compact Riemannian manifold $M$ through the application $T:G\times M\to M$.
For each $g\in G$, we define the transformation $T(g):M\to M$ by $T(g)x=T(g,x)$. For each $g\in G$, the transformation $T(g)$ is a diffeomorphism, $T(1)$ is the identity transformation on $M$ (where $1$ denotes the identity of $G$), and for $g,h\in G$, we have $T(g\cdot h)=T(g)\circ T(h)$ (denoting by $\cdot$ the group operation in $G$).

Now consider cocycles with values in a topological group $\Gamma$ (with an operation denoted by $+$, not necessarily commutative) with identity $e$. A \emph{cocycle $\alpha$ over $T$} with values in $\Gamma$ is a continuous transformation $\alpha:G\times M\to \Gamma$ such that \begin{equation*} \alpha \left(g_{2}g_{1},x\right) =\alpha \left( g_{2},T(g_{1})x\right) +\alpha \left(g_{1},x\right) \end{equation*} for each $x\in M$ and $g_{1},g_{2}\in G$.

Similarly, two cocycles $\alpha$ and $\beta$ over $T$ are said to be \emph{cohomologous} if there exists a continuous transformation $\Phi :M\to \Gamma$ such that \begin{equation} \alpha \left(g,x\right) =\Phi \left(T(g)x\right) +\beta \left( g,x\right) -\Phi (x), \label{quarentaum} \end{equation} for each $x\in M$ and $g\in G$, where $-$ represents the inverse in $\Gamma$.

A cocycle $\alpha :G\times M\to \Gamma$ is called a \emph{coboundary}, or \emph{cohomologically trivial}, if the equation \begin{equation} \alpha \left(g,x\right) =\Phi \left(T(g)x\right) -\Phi (x), \label{quarentadois} \end{equation} has a continuous solution $\Phi :M\to \Gamma$. This corresponds to the case where $\alpha$ is cohomolo-gous to the trivial cocycle $\beta (g,x)=e$. The transformation $\Phi$ is then called a \emph{trivialisation} of the cocycle $\alpha$. For a cocycle $\alpha$ to be cohomologically trivial, it must satisfy the identity
\begin{equation} 
\label{antonio} \alpha (g,x)=e \text{ for each } x\in M \text{ and } g\in G \text{ such that } T(g)x=x. 
\end{equation}

Note that for an Abelian group $\Gamma$, Equation \eqref{quarentaum} is equivalent to \begin{equation*} 
\alpha \left(g,x\right) -\beta \left(g,x\right) =\Phi \left( T(g)x\right) -\Phi (x). 
\end{equation*}

The study of cocycle cohomology has proven advantageous when considering groups with invariant metrics, the main examples being Lie groups. Let us now briefly recall some basic definitions from the theory of Lie groups.

A \emph{Lie group} $\Gamma$ is a differentiable manifold with a group structure such that, for every $g_{0}\in\Gamma$, the transformations $L_{g_{0}}:\Gamma\to \Gamma$ and $R_{g_{0}}:\Gamma\to \Gamma$, defined by $L_{g_{0}}(g)=g_{0}+g$ and $R_{g_{0}}(g)=g+g_{0}$, called \emph{left} and \emph{right translations}, respectively, are differentiable, as is the map $g\mapsto -g$. The left and right translations naturally induce, at each point $g\in\Gamma$, linear transformations $d_{g_{0}}L(g):T_{g}\Gamma\to T_{g_{0}+g}\Gamma$ and $d_{g_{0}}R(g):T_{g}\Gamma\to T_{g+g_{0}}\Gamma$ on tangent spaces. A vector field $X$ on $\Gamma$ is said to be \emph{left-invariant} if $d_{g_{0}}L(g)X(g)=X(g_{0}+g)$, and \emph{right-invariant} if $d_{g_{0}}R(g)X(g)=X(g+g_{0})$. The vector space $\mathcal{G}$ of left-invariant vector fields on $\Gamma$, equipped with a bilinear, anti-symmetric operation $[\cdot ,\cdot] :\mathcal{G}\times \mathcal{G}\rightarrow \mathcal{G}$ satisfying 
$$[X,[Y,Z]]+[Y,[Z,X]]+[Z,[X,Y]]=0,$$ 
is called the \emph{Lie algebra of $\Gamma$}. A Riemannian metric on $\Gamma$ is called \emph{left-invariant} if $L_{g_{0}}$ is an isometry for all $g_{0}\in\Gamma$, that is, $\left\langle g_{0}+g_{1},g_{0}+g_{2}\right\rangle =\left\langle g_{1},g_{2}\right\rangle$. A \emph{right-invariant metric} can similarly be defined by replacing $L_{g_{0}}$ with $R_{g_{0}}$. In any Lie group $\Gamma$, the choice of a symmetric, non-degenerate bilinear form $\left\langle \cdot,\cdot\right\rangle$ on the tangent space at the identity, $T_{e}\Gamma$, defines a left-invariant Riemannian metric on $\Gamma$ via \begin{equation*} 
\left\langle U,V \right\rangle_{g_{0}}=\left\langle d_{-g_{0}}L(g_{0})U,d_{-g_{0}}L(g_{0})V \right\rangle, 
\end{equation*} 
for all $U,V\in T_{g_{0}}\Gamma$. If the metric is also right-invariant, it is called \emph{bi-invariant}. Given $g\in\Gamma$, consider the conjugation map defined by $a\mapsto -g+a+g$ for all $a\in\Gamma$, and the transformation $Ad:\Gamma \to Aut(\mathcal{G})$ given by the derivative of conjugation by $g$, defined as $Ad(g)x=(d_{-g}R(g)\circ d_{g}L(e))x$ \, for $x\in\Gamma$.

The study of cohomology of cocycles over Anosov diffeomorphisms and flows (that is, hyperbolic actions of $G=\mathbb{Z}$ and $G=\mathbb{R}$, respectively) was initiated by Livschitz. In \cite{L1}, as discussed in \cite{RL2}, it was proven that a H\"{o}lder cocycle with values in $\Gamma =\mathbb{R}$ satisfying the conditions in \eqref{antonio} is cohomologous to the trivial cocycle via a H\"{o}lder solution $\Phi$. The same result was extended to $\Gamma$ as a connected Lie group admitting a bi-invariant metric.

The search for H\"{o}lder solutions to Equation \eqref{quarentaum} is discussed in \cite{P} and \cite{S3}.

Now, let us consider the measurable case. Let $T$ be an Anosov diffeomorphism, and $\mu$ a $T$-invariant measure equivalent to the Lebesgue measure. Livschitz showed in \cite{L2} that, for any H\"{o}lder cocycle $\alpha$ with real values, a measurable function $\Phi :M\to \mathbb{R}$ satisfying \eqref{quarentadois} almost everywhere is almost everywhere equal to a H\"{o}lder function $\widetilde{\Phi}$, for which $\alpha \left(g,x\right) =\widetilde{\Phi }\left(T(g)x\right) -\widetilde{\Phi }(x)$ everywhere. The same applies to Anosov flows. Pollicott and Walkden state in \cite{PW} that the proof presented by Livschitz extends to any connected Lie group with a bi-invariant metric. An analogous result for Equation \eqref{quarentaum} for compact Lie groups was proven by Parry and Pollicott in \cite{PP2}.

We now describe recent work by Pollicott and Walkden, which considers Lie groups that do not necessarily have a bi-invariant metric.

Let $f:M\to M$ be a $C^{1}$ diffeomorphism, and let $\Lambda$ be a locally maximal hyperbolic set for $f$. For a continuous function $g:\Lambda\to \Gamma$, consider the supremum 
\begin{equation*} 
\sup\left\{h_{m}(f)+ \int_{\Lambda}g,dm :m \text{ is an $f$-invariant probability measure}\right\}, 
\end{equation*} 
where $h_{m}(f)$ is the entropy of $f$ relative to $m$. If $g$ is H\"{o}lder, then this supremum is attained by some $f$-invariant probability measure, which we call the \emph{equilibrium measure of $g$} (see, for example, \cite{KH}). The hyperbolicity of $f\left|{\Lambda }\right.$ can be characterized in terms of the Mather spectrum of $f$. Let $\chi (\Lambda )$ denote the Banach space of continuous vector fields on $\Lambda$. Define the transformation $f^{\ast }:\chi (\Lambda )\to \chi (\Lambda )$ by
\begin{equation*} \left(f^{\ast }v\right)(x)=df(v(f^{-1}x)). 
\end{equation*} 
The Mather spectrum of $f$ is the spectrum of $f^{\ast }$ acting on the complexification of $\chi (\Lambda )$. If the Mather spectrum of $f\left|{\Lambda }\right.$ is contained in the set 
\begin{equation*}
\left\{ z\in \mathbb{C}:0<\left| z\right| <\lambda_{s}\right\} \cup \left\{ z\in \mathbb{C}:\lambda_{u}<\left| z\right| <\infty \right\}    
\end{equation*} 
for $\lambda_{s}<1<\lambda_{u}$, then we say that the function $\varphi :\Lambda \to \Gamma $ \emph{satisfies a partial hyperbolicity hypothesis} if we can choose $\lambda_{s}$ and $\lambda_{u}$ such that $\lambda_{s}<\mu_{s}\leq 1\leq \mu_{u}<\lambda_{u}$, where 
\begin{equation*} \mu_{s}=\limsup_{n\to \infty }\left(\sup_{x\in \Lambda }| Ad(\varphi (f^{n-1}x)\cdots\varphi (fx)\varphi (x))| \right)^{\frac{1}{n}}, 
\end{equation*} 
\begin{equation*} \mu_{u}=\liminf_{n\to \infty }\left(\sup_{x\in \Lambda }| Ad(\varphi (f^{n-1}x)\cdots\varphi (fx)\varphi (x))^{-1}| \right)^{-\frac{1}{n}}. 
\end{equation*}
Let be
\begin{equation*}
\widetilde{\theta}=\max \left\{ \frac{\log \mu_{s}}{\log \lambda_{s}},\frac{\log \lambda_{u}}{\log \mu_{u}}\right\} <1.
\end{equation*}

Assuming this partial hyperbolicity hypothesis, Pollicott and Walkden generalized in \cite{PW} Livschitz's H\"{o}lder regularity result to any connected Lie group (not necessarily possessing a bi-invariant metric).

\begin{theorem}
\label{teo10} Let $\Lambda$ be a compact locally maximal hyperbolic set of a $C^{1}$ diffeomorphism $f$, and $m$ an equilibrium measure of a H\"{o}lder function. Let $\Gamma$ be a connected Lie group and $\varphi :\Lambda \to \Gamma$ be a H\"{o}lder function with exponent $\theta \in (\widetilde{\theta},1)$ satisfying a partial hyperbolicity hypothesis. If $\Phi$ is a measurable solution of the cohomological equation $\varphi =\Phi \circ f-\Phi$ $m$-almost everywhere, then $\Phi$ is $m$-almost everywhere equal to a H\"{o}lder transformation $\widetilde{\Phi}$ for which $\varphi =\widetilde{\Phi} \circ f-\widetilde{\Phi}$ everywhere.
\end{theorem}

We note that this result does not guarantee the existence of solutions to the cohomological equation. We now describe a sufficiently general additional condition under which existence can be ensured. Let $\varphi :\Lambda \to \Gamma$ be a H\"{o}lder function satisfying a partial hyperbolicity hypothesis such that $\lambda_{s}<\mu_{u}^{-1}$ and $\mu_{s}^{-1}<\lambda_{u}$. We define the constant
\begin{equation*}
\widetilde{\theta^{\prime}}=\max \left\{ \widetilde{\theta},\frac{\log
\mu_{s}^{-1}}{\log \lambda_{u}},\frac{\log \mu_{u}}{\log
\lambda_{s}^{-1}}\right\} <1.
\end{equation*}
Still in \cite{PW}, Pollicott and Walkden proved the following result.

\begin{theorem}\label{teo11}
Let $\Lambda$ be a compact locally maximal hyperbolic set of a $C^{1}$ diffeomorphism $f$. Let $\Gamma$ be a connected Lie group and $\varphi :\Lambda \to \Gamma$ a H\"{o}lder function with exponent $\theta\in(\widetilde{\theta^{\prime}},1)$ satisfying a partial hyperbolicity hypothesis with $\lambda_{s}<\mu_{u}^{-1}$ and $\mu_{s}^{-1}<\lambda_{u}$. If
\begin{equation*}
\sum_{i=0}^{n-1}\varphi (f^{i}x)=0
\end{equation*}
whenever $f^{n}x=x$, then there exists a H\"{o}lder solution $\Phi :\Lambda \to \Gamma$ for the cohomological equation $\varphi =\Phi \circ f-\Phi$.
\end{theorem}

Pollicott and Walkden initially addressed the case where $\Gamma$ is a connected solvable Lie group (recall that a group is \emph{solvable} if there exist subgroups $\Gamma =\Gamma _{0}>\Gamma_{1}>\cdots>\Gamma_{n}=\left\{ e\right\}$ such that the quotients $\Gamma _{i}/ \Gamma _{i+1}$ are abelian, $0\leq i\leq n-1$) and then considered the more general case. Specifically, they showed that for solvable groups it is not necessary to assume any partial hyperbolicity hypothesis. To prove the results, cohomological equations over topological Markov chains and their suspensions  are considered (see \cite{RL1}), which allows solving first the equation in symbolic dynamics. It is worth noting that the work of Pollicott and Walkden is based on the articles \cite{PP1} and \cite{PP2}, where cocycles taking values in a compact Lie group over Anosov diffeomorphisms were studied. This problem was also the subject of several earlier articles preceding \cite{PW}, namely, \cite{NT3} by Ni\c{t}ic\u{a} and T\"{o}r\"{o}k, \cite{P} by Parry, and \cite{W} and \cite{W1} by Walkden.

Another generalization of the H\"{o}lder Livschitz's results was obtained by Ni\c{t}ic\u{a} and T\"{o}r\"{o}k for cocycles taking values in diffeomorphism groups of class $C^{k}$, $k=1,2,\ldots,\infty ,\omega$. Given a differentiable manifold $N$, a cocycle $\alpha :\mathbb{Z}\times M\to \text{Diff}^{k}(N)$ is said to be of class $C^{k}$ ($k=1,2,\ldots,\infty$) if the function $\varphi =\alpha (1,\cdot ):M\to \text{Diff}^{k}(N)$ is a function of class $C^{k}$. For $\Gamma =\text{Diff}^{k}(N)$, these authors proved in \cite{NT1} and \cite{NT2} that a cocycle with values close to the identity over a diffeomorphism $f$ with a hyperbolic set is cohomologous to the trivial cocycle provided that the conditions in \eqref{antonio} are satisfied. The proof strategy is analogous to the one used by Livschitz. The main difference is that in the group $\text{Diff}^{k}(N)$, the natural metric is neither left- nor right-invariant. However, in a sufficiently small neighborhood of the identity, it can be shown that this metric is ``almost invariant,'' which is sufficient for the proof.

\section{The Regularity Problem}

\qquad Another important question is what can be said about the regularity of the solution $\Phi$ of the cohomological equation when the cocycle $\alpha$ and the action $T$ are of class $C^{k}$, $k=1,2,\ldots,\infty,\omega$.

This question was also first studied by Livschitz. He showed in \cite{L1} that for $C^{1}$ cocycles $\alpha$, the solution $\Phi$ is still of class $C^{1}$. Livschitz also showed, using the decay of Fourier coefficients, that if $\alpha$ is a $C^{\infty }$ cocycle (respectively $C^{\omega }$) over certain linear actions on the torus, then the same holds for the solution $\Phi$ (see \cite{L2}). Subsequently, Guillemin and Kazhdan found $C^{\infty }$ solutions $\Phi$ for the case of $C^{\infty }$ geodesic flows on surfaces with negative curvature \cite{GK1,GK2}. In turn, Collet, Epstein, and Gallavotti proved in \cite{CEG} an analytic version also for geodesic flows but only on surfaces of constant negative curvature. All the aforementioned $C^{\infty }$ and $C^{\omega}$ regularity results were obtained through some form of generalized harmonic analysis, and as such, require a special structure on the manifold. Only in 1986, using a geometric argument, was a general result for $C^{\infty }$ regularity presented by de la Llave, Marco, and Moriy\'{o}n \cite{LMM}.

\begin{theorem}\label{teo12}
Let $M$ be a compact manifold, $f:M\to M$ a topologically transitive $C^{\infty }$ Anosov diffeomorphism, and $\varphi:M\to \mathbb{R}$ a $C^{\infty }$ function. Then, the following statements are equivalent:
\begin{enumerate}
\item There exists a $C^{\infty }$ function $\Phi :M\to \mathbb{R}$ satisfying
\begin{equation*}
\varphi =\Phi \circ f-\Phi;
\end{equation*}
\item For any periodic point $x$ of period $n$,
\begin{equation*}
\sum_{i=0}^{n-1}\varphi (f^{i}x)=0.
\end{equation*}
\end{enumerate}
Moreover, if $f$ is analytic, as well as the stable and unstable bundles and the function $\varphi$, then $\Phi$ is analytic.
\end{theorem}

These authors also established a corresponding result for flows.

\begin{theorem}
\label{teo13}
Let $M$ be a compact manifold, $\Psi =\left\{ \psi^{t}\right\}_{t\in \mathbb{R}}$ a topologically transitive $C^{\infty }$ Anosov flow, and $\varphi :M\to \mathbb{R}$ a $C^{\infty }$ function. Then, the following statements are equivalent:
\begin{enumerate}
\item There exists a $C^{\infty }$ function $\Phi :M\to \mathbb{R}$ satisfying
\begin{equation*}
\varphi(x)=\frac{d}{dt}\Phi (\psi^{t}x)|_{t=0};
\end{equation*}
\item For any periodic orbit $\left\{\psi^{t}x\right\}_{t\in \mathbb{R}}$ of period $T$, we have
\begin{equation*}
\int_{0}^{T}\varphi(\psi^{t}x)\,dt=0.
\end{equation*}
\end{enumerate}
Moreover, if $\Psi$ is analytic, as well as the stable and unstable bundles and the function $\varphi$, then $\Phi$ is analytic.
\end{theorem}

In the proof, an important general result of harmonic analysis is used, which states that if a function is of class $C^{\infty }$ along two absolutely continuous foliations of class $C^{\infty }$ with Jacobians having certain regularity properties, then the function is globally of class $C^{\infty }$. This result is proved in the article using the theory of elliptic operators and the absolute continuity of the Jacobian along the stable and unstable foliations. Since then, two new proofs for $C^{\infty }$ Anosov systems have emerged. One is due to Journ\'{e} in \cite{J}. An alternative approach was proposed by Hurder and Katok in \cite{HK}, based on an unpublished idea by Toll. Building on the approach presented in \cite{HK}, de la Llave later established the analytic case in 1997 \cite{Ll}. On the other hand, continuing the work developed by Ni\c{t}ic\u{a} and T\"{o}r\"{o}k in \cite{NT2,NT3}, Katok, Ni\c{t}ic\u{a}, and T\"{o}r\"{o}k extended in \cite{KNT} the regularity results to cocycles taking values in Lie groups and diffeomorphism groups.

\section{Cohomology in Higher Dimensions}

\qquad Another direction of research is the study of cohomology in higher dimensions for group actions such as $\mathbb{Z}^{k}$ and $\mathbb{R}^{k}$ for $k\geq 2$, as well as their "non-invertible" versions $\mathbb{Z}_{+}^{k}=\left\{(n_{1},\ldots,n_{k}):n_{1},\ldots,n_{k}\in \mathbb{N}\right\}$ and $\mathbb{R}_{+}^{k}=\left\{(x_{1},\ldots,x_{k}):x_{1},\ldots,x_{k}\in\mathbb{R}_{+}\right\}$.

Let $M$ be a compact manifold and let $T$ be an action of $\mathbb{Z}_{+}^{k}$ on $M$ generated by (not necessarily invertible) transformations $F_{1},\ldots,F_{k}:M\to M$ of class $C^{\infty }$ that commute with each other, that is, $F_{i}\circ F_{j}=F_{j}\circ F_{i}$ for all $i$ and $j$. For each $n\in \left\{1,\ldots,k\right\}$, an $n$\emph{-cochain} $T$ on $M$ with values in $\mathbb{R}^{l}$ $\left(l\geq 1\right)$ is any function
\begin{equation*}
\alpha :(\mathbb{Z}_{+}^{k})^{n}\times M\to \mathbb{R}^{l}
\end{equation*}
that is multi-linear and antisymmetric in the first $n$ variables and of class $C^{\infty }$ in the last variable. Since a multi-linear function is determined by its coefficients, the function $\alpha$ can be viewed as a $C^{\infty }$ function
\begin{equation*}
\alpha :M\to ( \mathbb{R}^{l})^{\binom{n}{k}}
\end{equation*}
with components indexed by $i_{1}<\cdots<i_{n}$, $i_{1},\ldots,i_{n}\in \left\{ 1,\ldots,k\right\}$.

Consider an operator $\mathcal{D}$, referred to as the \emph{coboundary operator}, which acts on each $n$-cochain $\alpha$ to produce an $(n+1)$-cochain $\mathcal{D}\alpha$ whose $\binom{k}{n+1}$ components are given by
\begin{equation*}
\left(\mathcal{D}\alpha \right)_{i_{1}\cdots i_{n+1}}\left(x\right) =\sum_{j=1}^{n+1}\left(-1\right)^{j+1}[
\alpha_{i_{1}\cdots\widehat{i_{j}}\cdots i_{n+1}}(F_{i_{j}}x)-
\alpha_{i_{1}\cdots\widehat{i_{j}}\cdots i_{n+1}}\left(x\right)] .
\end{equation*}
Here, the notation $\widehat{i_{j}}$ indicates that the index $i_{j}$ is omitted.

It is straightforward to verify that $\mathcal{D}^{2}=0$, enabling the introduction of a notion of cohomology associated with the coboundary operator. The cohomology of the $n$-cochain is referred to as the $n$\emph{-th cohomology} $C^{\infty }$ of the action $T$. The action $T$ of $\mathbb{Z}_{+}^{k}$ on $M$ naturally induces orbits in $M$. Let $P_{T}$ denote the set of all periodic orbits of the action $T$, that is, the finite orbits of $T$. For a periodic orbit $\mathcal{O}\in P_{T}$, let $\mu_{\mathcal{O}}$ be the unique normalized $T$-invariant measure associated with $\mathcal{O}$. Integrating $\alpha$ with respect to this measure produces a $n$\emph{-cocycle over} $T$ that is independent of the last variable, which determines an element $\left[ \alpha \right]_{\mathcal{O}}\in H^{n}\left(
\mathbb{Z}_{+}^{k};\mathbb{R}^{l}\right)$ in the $n$-th $C^{\infty }$ cohomology class with respect to the coboundary operator. Each $\left[ \alpha \right]_{\mathcal{O}}$ is a cohomological invariant of $\alpha$, meaning that $\left[ \alpha \right]_{\mathcal{O}}=\left[ \beta \right]_{\mathcal{O}}$ whenever
\begin{equation} \label{historia}
\alpha \left(t,x\right) =\beta \left(t,x\right) +\mathcal{D}\Phi
\left(t,x\right),
\end{equation}
since $[\mathcal{D}\Phi ]_{\mathcal{O}}=0$.

The action $T$ is said to satisfy the \emph{$C^{\infty}$ Livschitz's property for $n$-cocycles} if the set 
$$\left\{ \left[\alpha \right]_{\mathcal{O}}\mid \mathcal{O}\in P_{T}\right\}$$ 
constitutes a complete set of cohomology invariants for $n$-cocycles of class $C^{\infty }$. That is, if given $n$-cocycles of class $C^{\infty }$ $\alpha$ and $\beta$ such that $\left[ \alpha \right]_{\mathcal{O}}=\left[ \beta \right]_{\mathcal{O}}$ for every periodic orbit $\mathcal{O}$, there exists a $\left(n-1\right)$-cochain $\Phi$ that is a solution of Equation \eqref{historia}. In this case, the $n$-cocycles of class $C^{\infty }$ $\alpha $ and $\beta $ are said to be \emph{$C^{\infty }$-cohomologous}. An $n$-cocycle $\alpha \left(t,x\right) $ of class $C^{\infty }$ is called \emph{cohomologically trivial} if it is cohomologous to a constant $n$-cocycle $\beta (t)$. The $\left(n-1\right)$-cochain $\Phi $ that is a solution of Equation \eqref{historia} is called a \emph{trivialization} of $\alpha$. An $n$-cocycle $\alpha $ of class $C^{\infty }$ cohomologous to the trivial cocycle $\beta (t)=0$ for all $t\in \mathbb{Z}_{+}^{k}$ is called a \emph{coboundary}, and equation
\begin{equation*}
\alpha \left(t,x\right) =\mathcal{D}\Phi \left(t,x\right)
\end{equation*}
is called the \emph{cohomological equation}. Similarly, $C^{\infty }$ cohomology can be defined for an action $T$ of $G=\mathbb{Z}^{k}$ or $G=\mathbb{R}^{k}$. In the case $G=\mathbb{R}^{k}$, the $n$-cochains are vector fields of $n$-differential forms, the cocycles correspond to vector fields of closed forms, and the coboundary operators $\mathcal{D}$ are given by restrictions to the orbit foliation of $(n-1)$-differential forms of class $C^{\infty }$ globally defined.

In \cite{KK}, Katok and Katok used a version of Veech's method in \cite{V} to study $C^{\infty }$ cohomology for actions of hyperbolic automorphisms of the torus $\mathbb{T}^{N}$. They obtained, in particular, the following results.

\begin{theorem}
\label{teo15} Let $T$ be an action of $\mathbb{Z}^{k}$ by hyperbolic automorphisms of $\mathbb{T}^{N}$, and let $\alpha$ be a $k$-cocycle of class $C^{\infty }$ over $T$ with values in $\mathbb{R}^{l}$ $\left(l\geq 1\right)$ such that $\left[ \alpha \right]_{\mathcal{O}}=0$ for every $\mathcal{O}\in P_{T}$. Then, there exists a $\left(k-1\right)$-cochain $\Phi $ of class $C^{\infty }$ such that $\alpha =\mathcal{D}\Phi $.
\end{theorem}

\begin{theorem}
\label{teo16} Let $T$ be an action of $\mathbb{Z}^{k}$ by hyperbolic automorphisms of $\mathbb{T}^{N}$, and for $1\leq n\leq k-1$, let $\alpha $ be a $n$-cocycle of class $C^{\infty }$ over $T$ with values in $\mathbb{R}^{l}$ $\left(l\geq 1\right)$. Then, $\alpha $ is $C^{\infty }$-cohomologous to a constant cocycle $\beta$, that is,
\begin{equation*}
\alpha \left(t,x\right) =\beta \left(t\right) +\mathcal{D}\Phi
\left(t,x\right),
\end{equation*}
where $\Phi$ is a $\left(n-1\right)$-cochain of class $C^{\infty }$.
\end{theorem}
These results illustrate the power of the Livschitz's property for higher-order cohomology and complete the analysis of $C^{\infty }$ cohomology for actions of hyperbolic automorphisms.

We observe that for $n=k$ (Theorem~\ref{teo15}), the $C^{\infty}$ Livschitz's property holds: the \(k\)-th $C^{\infty}$ cohomology class of the action $T$ of $\mathbb{Z}^{k}$ is determined by the periodic orbits, meaning that periodic information is necessary and sufficient to identify coboundaries in the $k$-th $C^{\infty }$ cohomology. There is therefore a bijective correspondence between $k$-cocycles and functions on $M$ with values in $\mathbb{R}^{l}$. For the remaining possible values of \(n\) (Theorem~\ref{teo16}), $n\in \left\{1,\ldots,k-1\right\}$, each \(n\)-cocycle of class $C^{\infty }$ is cohomologous to a constant \(n\)-cocycle via an \((n-1)\)-cochain of class \(C^{\infty }\). Thus, the article \cite{KK} provides a complete description of \(C^{\infty }\) cohomology for cocycles with values in \(\Gamma =\mathbb{R}^{l}\) induced by an action of \(G=\mathbb{Z}^{k}\) through hyperbolic automorphisms of the toroidal manifold \(\mathbb{T}^{N}\). The main technique employed involves moving to a dual problem. The dual of a vector function on the torus \(\mathbb{T}^{N}\) is the collection of its Fourier coefficients, i.e., a vector function on \(\mathbb{Z}^{N}\). The dual cohomological equation is reduced to equations for the action on each orbit \(\mathcal{O}\) of the dual action. In particular, the differentiability of the original cocycle corresponds to super-polynomial decay of its Fourier coefficients.

For \(n=1\), Theorem~\ref{teo16} was established by Katok and Spatzier in \cite{KSp2} for the Anosov case and in \cite{KSp3} for the partially hyperbolic case. As noted by Katok and Katok, Livschitz's method is not valid for higher-order cohomology (Katok and Spatzier present in \cite{KSp2} an adaptation of Livschitz's method to establish the \(C^{\infty }\), \(C^{1}\), and H\"{o}lder Livschitz's properties but only for \(1\)-cocycles of Anosov actions of \(G=\mathbb{R}^{k}\)).

\section{Conclusion}

\qquad The research on cohomology for dynamical systems, particularly for cocycles and higher-dimensional actions, has made significant advances over the past decades. From Livschitz's foundational results to generalizations for \(C^{\infty }\) cocycles and actions on manifolds, this field continues to reveal profound connections between dynamical properties and cohomological invariants. Applications extend to hyperbolic dynamics, Lie groups, and beyond, providing a fertile ground for further exploration. Future research will likely uncover new techniques and frameworks to tackle even more complex dynamical and cohomological questions.

\end{document}